\documentclass[12pt]{article}
\usepackage[cp1251]{inputenc}
\usepackage[english ]{babel}
\usepackage{amssymb, amsmath}

\usepackage[dvips]{graphicx}
\begin{document}

\begin{center}
{\bf INTEGRAL GEOMETRY AND MIZEL'S PROBLEM}
\end{center}

\begin{center}
Yu.B.Zelinskii, M.V.Tkachuk, B.A.Klishchuk
\end{center}

\begin{center}
{\it Dedicated to professor Lech G\'{o}rniewicz in occasion of his
70 birthday}
\end{center}

\begin{center}
{\bf 1. Introduction}
\end{center}

A subject  treating in this article combines in one bundle some
questions of complex analysis, geometry and probability theory.
First investigations of geometric probability start from well
known Buffoons needle problem and related paradoxes.  Let a needle
be considered as a real line, and then the problem reduces to
finding some invariant measure of set relative to movement
(L. Santalo, G. Matheron, R. V. Ambarcumjan [1, 8, 12]).

Many questions of integral geometry are reduced to estimation of the
measure linear spaces crossing convex set. The finding of such
measure allows to the probabilistic estimations. Other problems
more close to geometry are estimations properties of set under
investigation if properties of its intersections with families of
some sets are well known:

1)  with planes of fixed dimension:

a) real case (G. Aumann, A. Kosi\'{n}ski, E. Shchepin [2, 7, 13]);

b) complex case (Yu. Zelinskii [18]);

2) with a set of vertices of an arbitrary  rectangle (A. Besicovitch, L. Danzer,
T. Zamfirescu, M. Tkachuk [3, 4, 17, 14]).

First problem is connected with the well known Ulam's problem from
Scottish book [9].

{\bf Ulam problem}. Let $M^n$ be $n$-dimensional manifold
and every section of $M^n$ by arbitrary hyperplane $L$ be homeomorphe
to $(n-1)$-dimensional sphere $S^{n-1}$. Is it true that $M^n$ is
$n$-dimensional sphere?

In the real case A.Kosi\'{n}ski solved this problem in 1962 [7].
L.Montejano received the repetition of this result in 1990 [10].
In complex case Yu.Zelinskii received similar result in 1993 [18].

The second problem is known in literature as Mizel's problem. Below
review of results  related to this problem will be given  and new unsolved
problems in this direction are discussed.

\begin{center}
{\bf 2. Mizel problem}
\end{center}

{\bf Mizel problem} (Characterization of a circle). A closed convex
curve such that, if three vertices of any rectangle lie on it, so
does the fourth, must be a circle.

In 1961, Besicovitch [3] solved
this problem. Later modified proof of this statement was presented
by Danzer [4], Watson [16], Koenen [6], Nash-Williams [11].

In 1989, T. Zamfirescu [17] proved the similar result for Jordan
curve (not convex a priory) and for rectangle with infinitesimal
relation between sides.
$$
\left| \frac{ac}{ab} \right| \leq \varepsilon > 0.
$$

In 2006, M. Tkachuk [14] received the most general result in
this area for an arbitrary compact set $C \subset \mathbb{R}^2$, where the
complement $\mathbb{R}^2 \backslash C$ is not connected.

Obviously that requirement to compactness is necessary, otherwise
straight line and some other sets, noted below, will satisfy the
assessed requirement. But if not to require partition to planes,
the following ensembles  will satisfy, for instance, that condition of the rectangle  : ensemble from three points of plane such that
triangle with vertices in these points will not be rectangular, a proper
arc of semicircle, set points of plane with rational (irrational)
coordinates.

Zamfirescu in his article [17] invited his readers to verify whether
indeed an arbitrary convex curve $ \Gamma $ of constant width $ d $
satisfying the infinitesimal rectangular property is a circle. In
this paper we shall give an affirmative answer to this question in
the following theorem:

{\bf Theorem 1.} The convex curve of constant width satisfying the
infinitesimal rectangular condition  is a circle.

Combining this result with the result of Zamfirescu [17] we obtain:

{\bf Theorem 2.} Every Jordan curve satisfying infinitesimal
rectangle condition is a circle.

We are only interested rectangles with the diagonal length $ d $, so
the infinitesimal condition may be replaced by the requirement that
the smaller rectangle side has length less than some $ \varepsilon
$.

In [14] it was proved that the curve $ \Gamma $ has a continuous
tangent.

We introduce the following notation: at each point $ x \in \Gamma $
denote by $ \partial \Delta_x $ the circle having with the curve $
\Gamma $ a common tangent at the point $ x $; in some neighborhood
of a point of the curve $\Gamma$ we assume that upward direction is
the direction along the the inner normal and according with this we
consider the right direction and  the left direction along the curve
$\Gamma$ ; $U_{\varepsilon} (x,\Gamma)$ is the
$\varepsilon$-neighborhood of $x$ on the curve $\Gamma$,
$U_{\varepsilon}^l (x,\Gamma)$ is the left
$\varepsilon$-neighborhood of $x$ on  the curve $\Gamma$, i.e. the
subset of $U_{\varepsilon} (x,\Gamma)$ each point of which lies to
the left of $x$; $U_{\varepsilon}^r (x,\Gamma)$ is the right
$\varepsilon$-neighborhood of $x$ on $\Gamma$.

Having in the proof provided by Besicovitch [3] a local character all
four lemmas are valid for the infinitesimal rectangle property. But
the result following from these lemmas changes: if at some point $ x
\in \Gamma $  the circle $ \partial \Delta_x $ intersects $ \Gamma $
at $ y $ and the distance between $x$ and $y$ is less than
$\varepsilon$, then in some neighborhood of $ x $ (and also of the
opposite point $ x ^ * $) the curve $ \Gamma $ is an arc of the
circle $ \partial \Delta_x $.

Consequently all points of the curve $ \Gamma $ can be divided into
five disjoint sets:

$ A = \{x \in \Gamma | U_{\varepsilon} (x,\Gamma) \backslash \{ x \}
\subset \Delta_x \}$ is the set of points of the curve $ \Gamma $,
such that the curve $ \Gamma $ lies in the open circular disk $
\Delta_x $ in $ \varepsilon $-neighborhood of $x \in A$ except the
point $ x $ itself;

$ B = \{x \in \Gamma | U_{\varepsilon} (x,\Gamma) \backslash \{ x \}
\subset \mathbb{R}^2 \backslash \overline{\Delta_x} \}$ is the set
of points of the curve $ \Gamma $ such that in their
$\varepsilon-$neighborhood the curve $\Gamma$ is situated outside of
the closed circular disk $\overline{\triangle_{x}}$ except the point
$x$ itself;

$ AB = \{x \in \Gamma | U_{\varepsilon}^l (x,\Gamma) \subset
\Delta_x, U_{\varepsilon}^r (x,\Gamma) \subset \mathbb{R}^2
\backslash \overline{\Delta_x} \}$ is the set of points of the curve
$ \Gamma $ such that in their $\varepsilon-$neighborhood the curve
$\Gamma$ lies in the open circular disk $\triangle_{x}$ to the left
of $x$ and outside of the closed circular disk
$\overline{\triangle_{x}}$ to the right of $x$ except the point $x$
itself;

$ BA = \{x \in \Gamma | U_{\varepsilon}^l (x,\Gamma) \subset
\mathbb{R}^2 \backslash \overline{\Delta_x} , U_{\varepsilon}^r
(x,\Gamma) \subset \Delta_x \}$ is the set of points of the curve $
\Gamma $ such that in their $\varepsilon-$neighborhood the curve
$\Gamma$ lies in the open circular disk $\triangle_{x}$ to the right
of $x$ and outside of the closed circular disk
$\overline{\triangle_{x}}$ to the left of $x$ except the point $x$
itself;

$ C $ is the set of points of the curve $ \Gamma $ each of which has
the neighborhood where curve coincides with an arc of the circle $
\partial \Delta_x $.

Hence $ \Gamma = A \bigcup B \bigcup AB \bigcup BA \bigcup C $.

Suppose that $ x_n \rightarrow x $, $ x_n \in \Gamma $ and let all $
x_n $ be in $ \varepsilon $-neighborhood of $ x $ to the left of it.
Assume that $ x_n \in A \bigcup BA $ and $x\in AB$ then in some
neighborhood of $ x $  the curve $ \Gamma $ lies above each circle $
\partial \Delta_ {x_n} $. $\Gamma$ is inside of the
circular disk $\overline{\triangle_{x}}$. This fact follows from the
continuity of the tangent to $\Gamma$. We obtain a contradiction
with $x\in AB$. So in some neighborhood of $ x $ to the left of it
there are no points of the sets $ A $ and $ BA $. In this
neighborhood we consider the similar sequence of points $ x_n \in B
\bigcup AB $ and as a consequence we find that in some neighborhood
of $ x $  the curve $ \Gamma $ lies under every circle $
\partial \Delta_ {x_n} $ and therefore outside of the circular
disk  $ \Delta_x $ which is impossible. So in some left
half-neighborhood of $ x $ our curve is an arc of a circle $
\partial \Delta_x $ which is also not compatible with the condition
$ x \in AB $. We conclude that $ AB = \emptyset $. Similarly $ BA =
\emptyset $.

Thus $ \Gamma = A \bigcup B \bigcup C $. Suppose that $ C =
\emptyset $, $ \Gamma = A \bigcup B $. Consider the sequence of
points $ x_n \in A $, $ x_n \rightarrow x $. If $ x \in B $ then as
we know the contradiction is obtained. Hence $x$ does belong to $A$
and $A$ is a closed set. Similarly $ B $ is closed. The curve
$\Gamma $ is connected and therefore either $ A $ or $ B $ is empty.
It means that $ \Gamma = A $ or $ \Gamma = B $ which implies that
the curve $ \Gamma $ either contains a circular disk $ \Delta_x $ or
is contained in a circular disk $ \Delta_x $ and it contradicts the
fact that $ \Gamma $ is the convex curve of constant width and has a
length $ \pi d $.

Hence $ C \neq \emptyset $ and there exists a point $ x \in C $.
Some neighborhood of $ x $ is also contained in the set $ C $. We
shall move along the curve $ \Gamma $ to the left from the point $ x
$ to the first point $y$ that does not belong to $C$. It is
obviously that $y$ cannot belong neither to $A$ nor to $B$ and
therefore $\Gamma = C$. Then by applying the Heine-Borel lemma we
conclude that $\Gamma$ is a circle and theorem 1 is proved.

The row of similar opened problems in the plane and in $n$-dimensional
case appears in connection with Mizel's problem.

{\bf Problem 1.} Let $C$ be closed Jordan curve in $\mathbb{R}^2$ and for arbitrary
algebraic closed curve $L$ of order $n$ from property that  {intersection
$C \cap L$ contains $m$ points} follows, that $C \cap L$ contains no less then $m+1$
points. Does there exists a number $m$, that from property above
follows that $C$ be algebraic curve of order $n$?

{\bf Problem 2.} Let in previous question $L$ be a circle and $m = 3$, is it true,
that $C$ also be a circle?

{\bf Problem 3.} Let in problem 1 $L$ be an ellipse and $m = 4$, is it true,
that $C$ also be an ellipse?

{\bf Problem 4.} Will be a compact $C$ a sphere in $\mathbb{R}^n$, if $C$ divides the
space, and if from the belonging $n+1$ tops of the arbitrary
rectangular parallelepiped to a compact $C$, it follows that one
more top lies in $C$ too?

Last question is interesting even if $C$ be
$(n - 1)$-dimensional manifold or boundary of a convex set.

{\bf Problem 5.} Let $C$ be a $(n-1)$-dimensional manifold (or boundary of a
convex domain) in $\mathbb{R}^n$ and not exist $(n-1)$-dimensional sphere $S^{n-1}$
that intersection $C \cap S^{n-1}$ contains $n+1$ points exactly. Is it
true, that $C$ be a $(n-1)$-dimensional sphere?

\textbf{Corollary.} On two-dimensional plane $\mathbb{R}^2$ there
exist compact sets, which divide the plane and such that no circles
has with them in accuracy three crosses point. In particular class
of similar sets includes the Shottka sets and the Sierpi\'{n}ski's
carpet.

{\bf Problem 6.} Remain or not result of [14] true, if we will consider compact
 set $C \subset \mathbb{R}^2$, where the complement $\mathbb{R}^2 \backslash C$ is not connected?

{\bf Problem 7.} Are cited results and problems 1-6 true, if we will consider that
one point (top) on $C$ is fixed?

Next examples will show that in problems 1-3 it is impossible
instead of curve, in analogy with Tkachuk's result, consider compact set dividing plane.

{\bf Example.} We shall consider the domain $D$ on plane, bounded by circle $S^1$.
The Surge from it thick in the domain $D$ infinite ensemble of opened balls $D_i$,
which are not intersect pair wise even on border and also not intersect the circle $S^1$.

\begin{figure}[h]
\center{\includegraphics[width=0.4\linewidth]{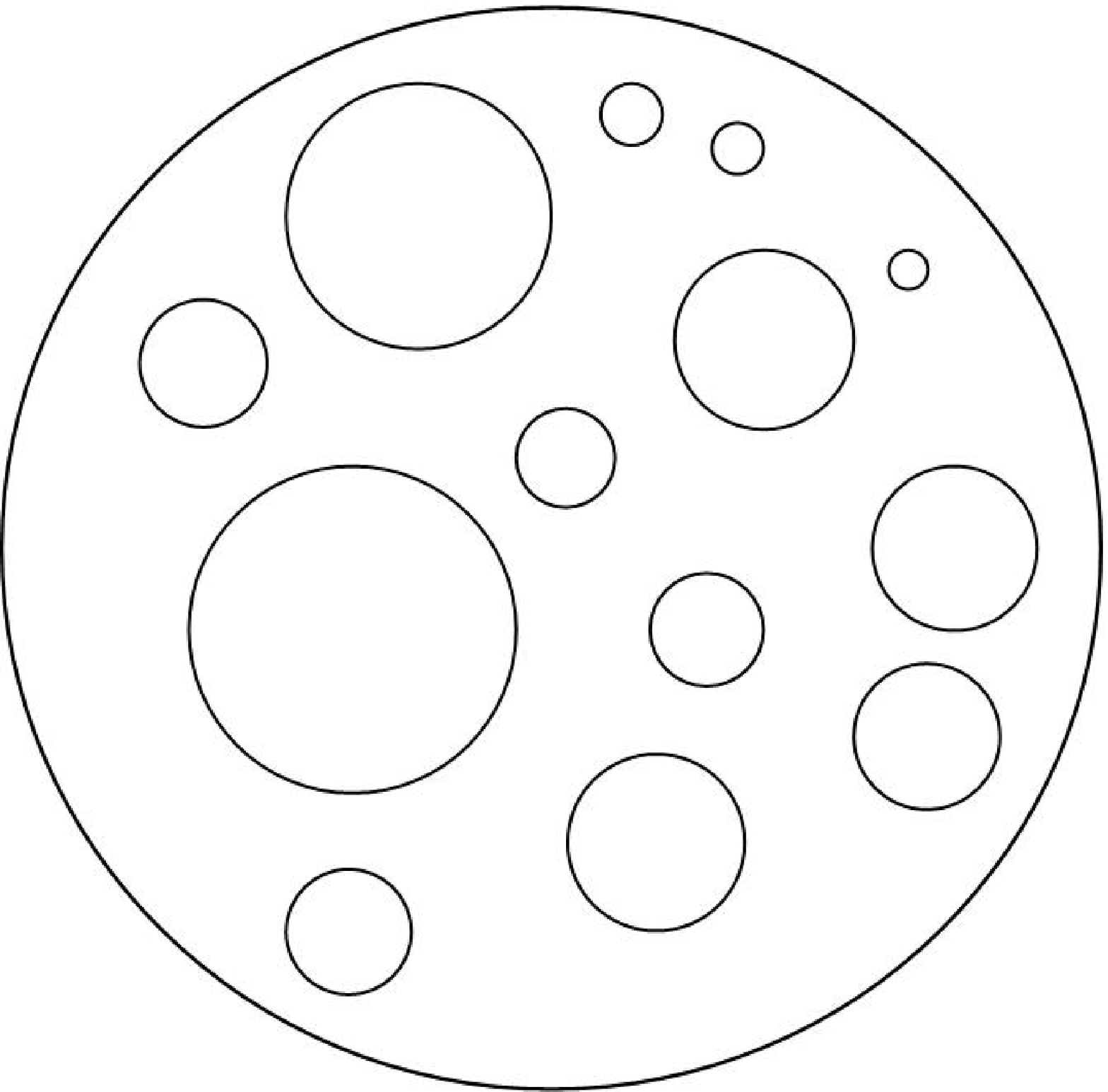} \\ Figure 1}
\label{ris:image}
\end{figure}

Then we receive fractal compact set $K=\bar{D}\backslash \cup D_i$
without interior points, which divides plane on countable set of
components. But it is easy to see that arbitrary circle can
intersect $K$ on only one point or on infinite number of points
(see Figure 1).

\begin{figure}[h!]
\center{\includegraphics[width=0.5\linewidth]{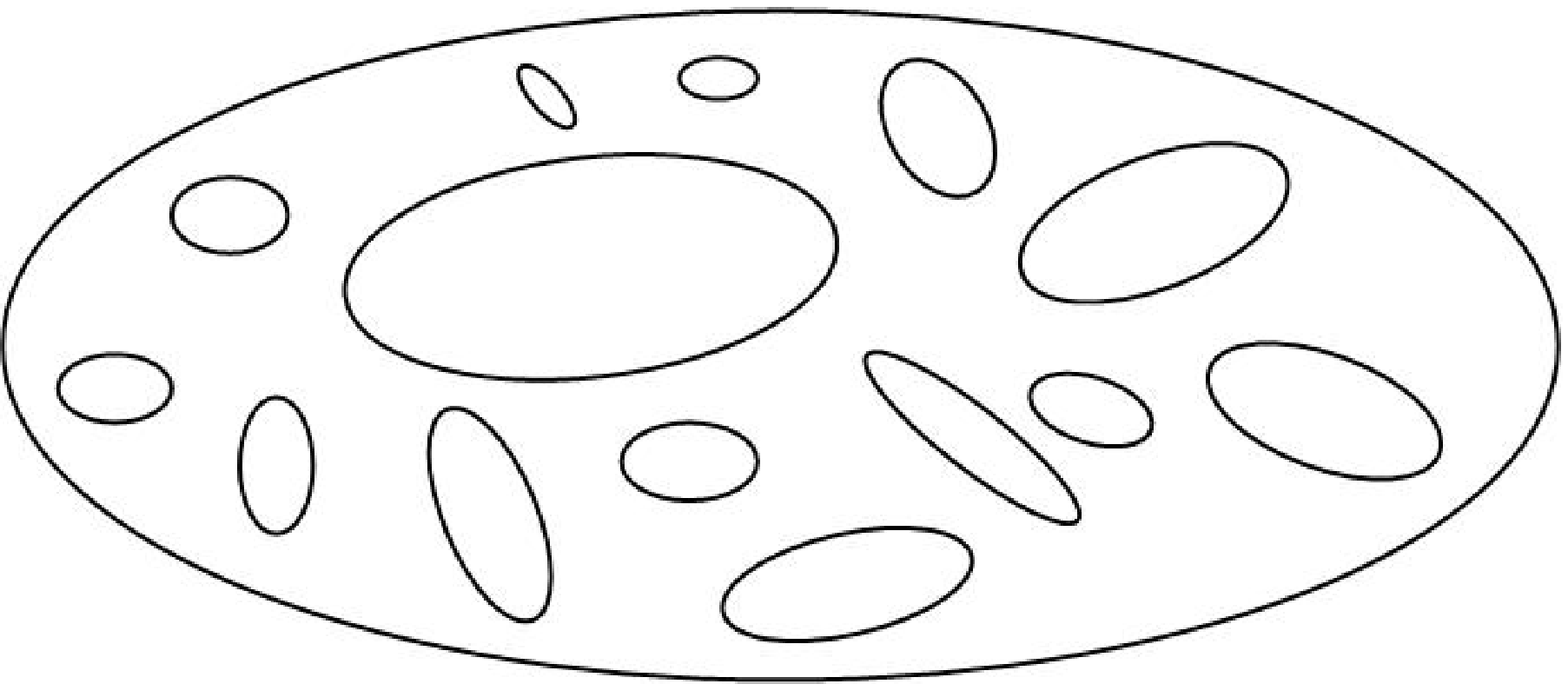} \\ Figure 2}
\label{ris:image}
\end{figure}

\begin{figure}[h!]
\center{\includegraphics[width=0.4\linewidth]{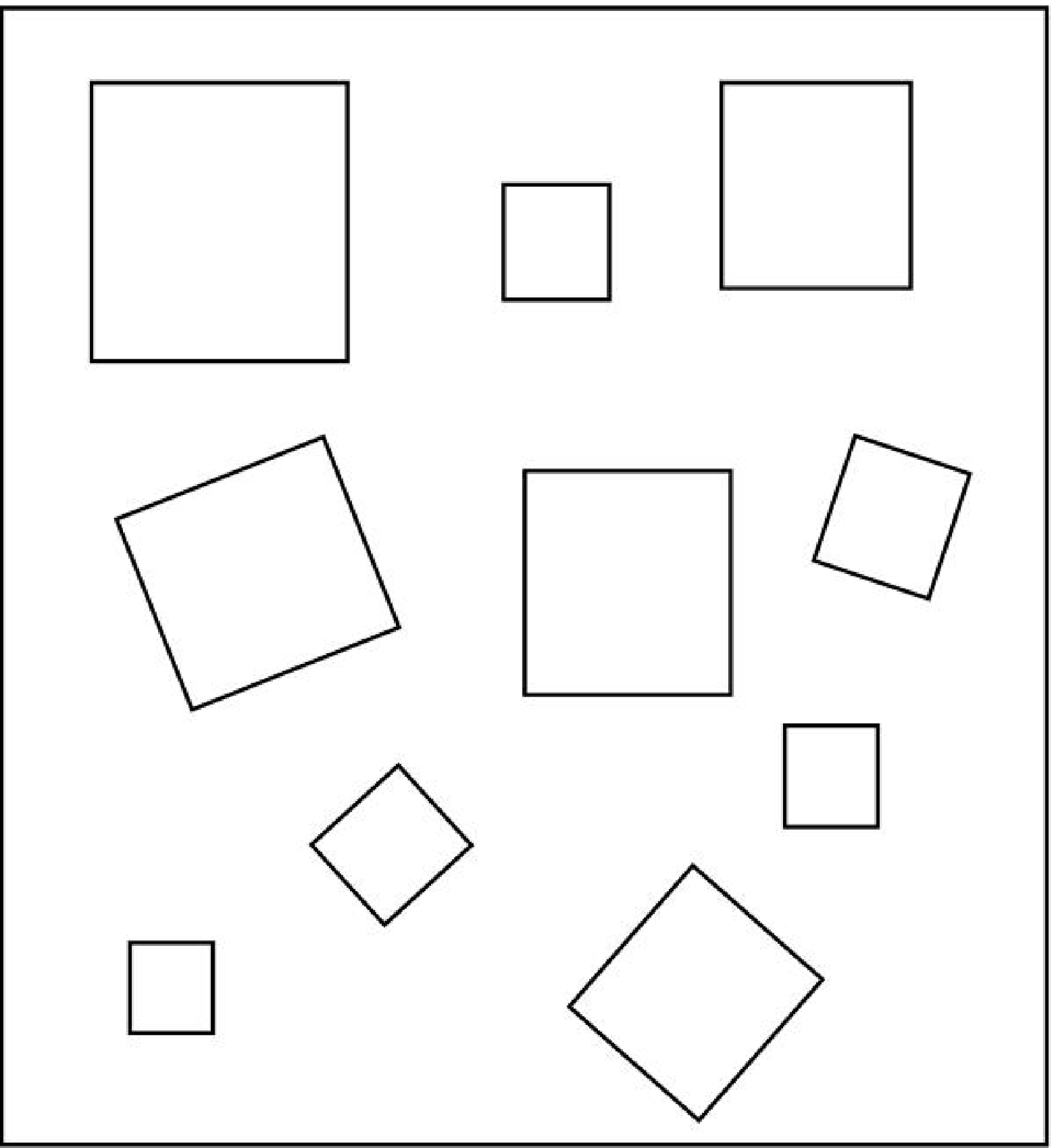} \\ Figure 3}
\label{ris:image}
\end{figure}

Other examples we will receive if domains $D (D_i)$ be domains
bounded by ellipses or squares. In this case intersections with
circle or ellipse can contain one, two, four or infinite many points
but never three points (see Figures 2 and 3). Those examples give
negative answer to problem 8 from [19].

Other close problems possible to find in the work of Gr\"{u}nbaum [5].

This investigation was partially supported by Tubitek-NASU grant
number 110T558

\end{document}